\documentclass[a4paper,11pt]{article}
\usepackage{amsmath} \usepackage{amsfonts} \usepackage{amssymb}\usepackage{latexsym}
\usepackage[english]{babel}
\usepackage{amscd}
\usepackage[dvips,final]{graphics}
\usepackage{fancyhdr}
\usepackage{locengli,math}
\usepackage{graphicx,pst-plot,psfig,pstricks,pst-node,pst-text,pst-tree,pst-char,pst-coil}

\newcommand{\Aut}{\textrm{Aut}}

\begin{document}

\begin{center}
\textbf{\LARGE{\textsf{On representations of braid groups determined by
directed graphs}}}
\footnote{
\textit{A.M.S. classification 2000: 16W30; 05C20; 05C90. }
}

\vskip1cm
\parbox[t]{14cm}{\large{
Philippe {\sc Leroux}}\\
\vskip4mm
{\footnotesize
\baselineskip=5mm
Institut de Recherche
Math\'ematique, Universit\'e de Rennes I and CNRS UMR 6625\\
Campus de Beaulieu, 35042 Rennes Cedex, France, pleroux@univ-rennes1.fr}}
\end{center}

\vskip1cm
{\small
\vskip1cm
\baselineskip=5mm
\noindent
{\bf Abstract:}
We prove that any non $L$-cocommutative finite Markov $L$-coalgebra yields at least two solutions
of the Yang-Baxter equation and therefore at least two representations
of the braid groups.
We conclude by a generalisation of these constructions to any coalgebra.
\section{Introduction}
In this article, $k$ is either the real field or the complex field.
Moreover, all the involved vector spaces will have a finite or
a denumerable basis.

The first part recalls the main notions on $L$-coalgebras introduced in
\cite{Coa} and developed in \cite{Coa}\cite{codialg1}\cite{perorb1}.
The second part uses these objects to construct representations of braid groups.
\section{$L$-coalgebras}
\begin{defi}{[$L$-coalgebra]}
A \textit{$L$-coalgebra} $G$ 
over a field $k$ is a $k$-vector space equipped with a right coproduct,
$\Delta: G \xrightarrow{} G^{\otimes 2}$ and a left coproduct, $\tilde{\Delta}: G \xrightarrow{} G^{\otimes 2}$, verifying the
coassociativity breaking equation 
$(\tilde{\Delta} \otimes id)\Delta = (id \otimes \Delta)\tilde{\Delta}$.
If $\Delta = \tilde{\Delta}$, the coalgebra is said  \textit{degenerate}. A $L$-coalgebra may have up to two counits, the right counit  
$\epsilon: G \xrightarrow{} k$, verifying $ (id \otimes \epsilon)\Delta = id$
and the left counit $\tilde{\epsilon}: G \xrightarrow{} k$, verifying $ ( \tilde{\epsilon} \otimes id)\tilde{\Delta} = id. $
\end{defi}
\begin{defi}{[finite Markov $L$-coalgebra]}
A {\it{finite Markov $L$-coalgebra}} $G$ with dimension $\dim G$
and with a basis $(v_i)_{1 \leq i \leq \dim G}$ is a $L$-coalgebra such that for all $v_i$ with $1 \leq i \leq \dim G$, $\Delta v_i = \sum_{k: \ v_i \otimes v_k \in I_{v_i}} \ w_{v_i}(v_i \otimes v_k) \ v_i \otimes v_k $
and $ \tilde{\Delta}v_i =\sum_{j: \ v_j \otimes v_i \in J_{v_i}} \ \tilde{w}_{v_i}(v_j \otimes v_i) \ v_j \otimes v_i $, where 
$I_{v_i},J_{v_i}$ are finite sets and $w_{v_i}: I_{v_i} \xrightarrow{} k$ and
$\tilde{w}_{v_i}: J_{v_i} \xrightarrow{} k$
are linear mappings called \textit{weight}.
\end{defi}
Let $V$ be a $k$-vector space.
Denote by $\tau$, the \textit{transposition} mapping, i.e. $V^{ \otimes 2} \xrightarrow{\tau} V^{ \otimes 2}$ such that $\tau(x \otimes y) = y \otimes x$ for all $x,y \in V$.
\begin{defi}{[$L$-cocommutativity]}
Let $G$ be a $L$-coalgebra. Denote its coproducts by $\tilde{\Delta}$ and $\Delta$. $G$ is said {\it{$L$-cocommutative}} iff for all $v \in G$, 
$(\Delta - \tau\tilde{\Delta})v = 0$.
\end{defi}
\begin{defi}{[Directed graph]}
A \textit{directed graph} $G$ is a quadruple \cite{Petritis}, $(G_{0},G_{1},s,t)$
where $G_{0}$ and $G_1$ are two denumerable sets respectively called the \textit{vertex set} and the \textit{arrow set}.
The two mappings, $s, \ t: G_1 \xrightarrow{} G_0$ are respectively called  \textit{source} and \textit{terminus}.
A vertex $v \in G_0$ is a \textit{source} (resp. a \textit{sink}) if $t^{-1}(\{v \})$ (resp. $s^{-1}(\{v\})$)
is empty. A graph $G$ is said {\it{locally finite}}, (resp. {\it{row-finite}}) if 
$t^{-1}(\{v\})$ is finite (resp. $s^{-1}(\{v\})$ is finite).
Let us fix a vertex $v \in G_0$.
Define the set $F_{v} :=\{a \in G_{1}, \ s(a)=v \}$. A \textit{weight} associated with
the vertex $v$ is a mapping $w_v: F_{v} \xrightarrow{} k$.
A directed graph equipped with a family of weights $w := (w_v)_{v \in G_0}$ 
is called a weighted graph.
\end{defi}
In the sequel, directed graphs will be supposed locally finite and row finite without sink and source.

We recall a theorem from \cite{Coa}.
Let $G$ be a directed graph equipped with a family of weights $(w_v)_{v \in G_0}$. Let us consider the free vector space generated by $G_0$.
The set $G_1$ is then viewed as a sub-vector space of $G_0^{\otimes 2}$. The 
mappings source and terminus are then linear mappings still called 
source and terminus $s,t: \ G_0^{\otimes2} \xrightarrow{} G_0$, such that $s(x \otimes y) = x$ and $t(x \otimes y) = y$ for all $x,y \in G_0$. The family of weights $(w_v)_{v \in G_0}$ is then viewed as 
a family of linear mappings from $F_v$ to $k$.
Let $v \in G_0$ and define the right coproduct $\Delta$ such that
$\Delta(v) := \sum_{i: a_i \in F_v} \ w_v(a_i) \ v \otimes t(a_i)$ and the left coproduct
$\tilde{\Delta}$ such that $\tilde{\Delta}(v) := \sum_{i: a_i \in P_v} \ w_{s(a_i)}(a_i) \ s(a_i) \otimes v$, where $P_{v} $ is the set $\{a \in G_{1}, \ t(a)=v \}$. Define, for all $v \in G_0$, the linear
mappings $\tilde{w}_v : P_v \xrightarrow{} k$
such that $\tilde{w}_v (a_i) = w_{s(a_i)}(a_i)$ for all $a_i \in P_v$.
With these definitions the vector space $G_0$ is a finite Markov $L$-coalgebra. 
In the sequel, $G_0$ will be identified with $G$. 
\NB \textbf{[Geometric representation]}
To yield a geometric support for coalgebras, we associate
with each tensor product $\lambda x \otimes y$, where $\lambda \in k$, appearing in the definition of the coproducts a directed arrow
$x \xrightarrow{\lambda}y$. The directed graph so obtained is then called the \textit{geometric support} of the coalgebra. The advantage of this
formalism is to generalise the notion of directed graph.
That is why we identify notions from directed graphs with notions from
$L$-coalgebras. 
We draw attention to the fact that a directed graph can be the
geometric support of different $L$-coalgebras.
\begin{exam}{}
The directed graph:
\begin{center}
\includegraphics*[width=4cm]{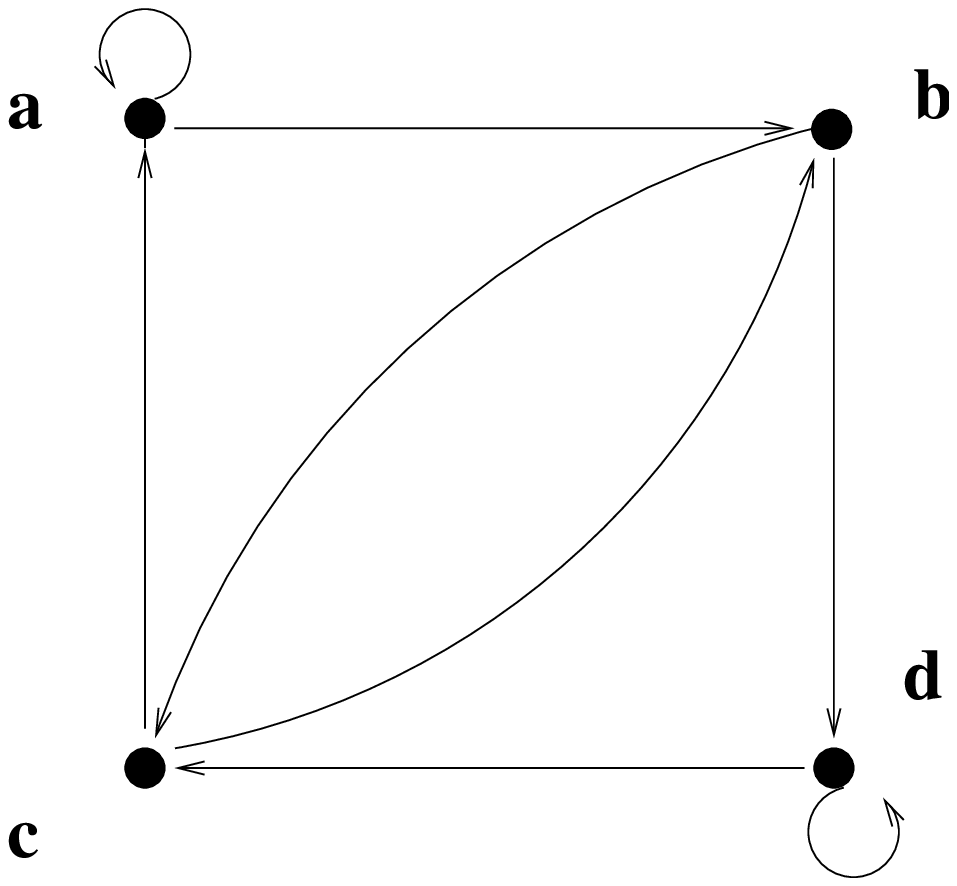}
\end{center}
\end{exam}
is the geometric support of the degenerate or coassociative $L$-coalgebra, generated by $a,b,c$ and $d$ and
described by the following coproduct:
$
\Delta a = a \otimes a + b \otimes c, \ \
\Delta b = a \otimes b + b \otimes d, \ \
\Delta c = d \otimes c + c \otimes a, \ \
\Delta d = d \otimes d + c \otimes b
$
and the geometric support of the finite Markov $L$-coalgebra, generated by $a,b,c$ and $d$ and described by the right coproduct:
$
\Delta_M a = a \otimes (a + b), \ \
\Delta_M b = b \otimes (c + d), \ \
\Delta_M c = c \otimes  (a +b), \ \
\Delta_M d = d \otimes (c + d)
$
and the left coproduct: 
$
\tilde{\Delta}_M a = (a+c) \otimes a,  \ \
\tilde{\Delta}_M b = (a +c) \otimes b , \ \
\tilde{\Delta}_M c = (b + d) \otimes c , \ \
\tilde{\Delta}_M d = (b+ d) \otimes d.
$
\NB
Let $G$ be a finite Markov $L$-coalgebra. If the family of weights
used for describing right and left
coproducts take values into $\mathbb{R}_+$ and
if the right counit $\epsilon: v \mapsto 1$ exists, then the geometric
support associated with $G$
is a directed graph equipped with a family of probability vectors. 
\section{Representation of braid groups determined by directed graphs}
\begin{defi}{[Yang-Baxter equation]}
Let us consider a $k$-vector space $V$. Let  $\hat{\Psi}$ be
an automorphism on $V^{\otimes 2}$, $\hat{\Psi}$ verifies
the \textit{Yang-Baxter equation} (YBE) if  \cite{Kassel}:
$$ (\hat{\Psi} \otimes id)(id \otimes \hat{\Psi})(\hat{\Psi} \otimes id) =(id \otimes \hat{\Psi})(\hat{\Psi} \otimes id)(id \otimes \hat{\Psi}).$$
Such a solution is also called a $R$-matrix. Let us denote by $S$ the set of solutions of YBE and by $\Aut (V)$ the linear automorphisms group of $V$.
\end{defi}
\NB
Let us recall that any solution of YBE supplies a representation of
braid groups, see for instance \cite{Kassel}\cite{Majid}.
The aim of this article is to show that markovian coproducts, used
to code the paths of directed graphs,
yield solutions of YBE and thus representation of braid groups. 
\begin{theo}
\label{tttqye}
Let $V$ be a $k$ vector space. The mapping $\ \widehat \ : \Aut (V) \xrightarrow{} S \times S$ defined by
$\Psi \mapsto (\hat{\Psi}_1, \ \hat{\Psi}_2)$ where
$\hat{\Psi}_1 := \tau(id \otimes \Psi)$ et $\hat{\Psi}_2 := \tau(\Psi  \otimes id )$ is injective.
\end{theo}
\Proof
The injectivity is trivial.
Let $\Psi$ be an automorphism on $V$ and $x,y,z \in V$.
Let us prove that $\hat{\Psi}_1$ is a solution of YBE.
\begin{eqnarray*}
 x \otimes y \otimes z &\xrightarrow{\hat{\Psi}_1 \otimes id}& \Psi(y) \otimes x  \otimes z \xrightarrow{ id \otimes \hat{\Psi}_1 }
\Psi(y) \otimes \Psi(z)  \otimes x \xrightarrow{\hat{\Psi}_1 \otimes id} \Psi^2(z) \otimes \Psi(y)  \otimes x. \\
x \otimes y \otimes z & \xrightarrow{id \otimes \hat{\Psi}_1}& x \otimes \Psi(z)  \otimes y \xrightarrow{ \hat{\Psi}_1 \otimes id}
\Psi^2(z) \otimes x  \otimes y \xrightarrow{id \otimes \hat{\Psi}_1} \Psi^2(z) \otimes \Psi(y)  \otimes x.
\end{eqnarray*}
$\hat{\Psi}_1$ is thus a solution of YBE. It is the same proof for  $\hat{\Psi}_2$. 
\eproof
\begin{theo}
Let $V$ be a $k$-vector space. If $A$ and $B$ are automorphisms on $V$
we define $c(A,B) := \tau(A \otimes B)$.
Let $\Psi_1$, $\Psi_2$, $\Psi_3$ and $\Psi_4$ be four automorphisms on $V$. 

\noindent
$c(\Psi_1, \Psi_2)$
is a solution of YBE iff  $\Psi_1 \Psi_2 =\Psi_2 \Psi_1$. If $c(\Psi_1, \Psi_2)$ and  $c(\Psi_3, \Psi_4)$ are
solutions of YBE and $[\Psi_1 \Psi_4, \Psi_2 \Psi_3] =0$  then
$\tau c(\Psi_1, \Psi_2) c(\Psi_3, \Psi_4)$ is still a solution of YBE.
\end{theo}
\Proof
Let $\Psi_1$, $\Psi_2$, $\Psi_3$ and $\Psi_4$ be four automorphisms on $V$.
The first claim is straightforward and so is the second one by noticing that
$\tau c(\Psi_1, \Psi_2) c(\Psi_3, \Psi_4) = c(\Psi_2 \Psi_3, \Psi_1 \Psi_4)$.
\eproof
\begin{defi}{[The companion graph]}
Let $G$ be a finite Markov $L$-coalgebra, equipped with a basis $(v_i)_{i = 1, \ldots, \dim G}$ and equipped with two coproducts $\Delta$ and $\tilde{\Delta}$.
These coproducts define
a directed graph $G$ without sink and source, where the basis vectors are the vertex of $G$, equipped with two family of weights $(w_{v_i})_{v_i \in G}$ and $(\tilde{w}_{v_i})_{v_i \in G}$. 
By definition, for all $v_i \in G$, there exists finite sets $I_{v_i}$ and $J_{v_i}$ such that
$\tilde{\Delta}v_i = \sum_{k: \ v_{k} \otimes v_i \in J_{v_i}} \tilde{w}_{v_{i}}(v_{k} \otimes v_i) \ v_{k} \otimes v_i$
and $\Delta v_i = \sum_{k: \ v_i \otimes v_{k} \in I_{v_i}} w_{v_i}(v_i \otimes v_{k}) \
v_i \otimes v_{k}$, for all $i$, with $1 \leq i \leq \dim G$.

\noindent
Let $H$ be a $k$-vector space such that dim $H$ = dim $G$, equipped with a 
basis
$(h_i)_{i = 1, \ldots, \dim G}$. 
With each vector $v_i \in G, \ 1 \leq i \leq \dim G$, let us associate an unique $h_i \in H$. Denote by
$G_* := G \oplus H$, the associated finite Markov
$L$-coalgebra, such that $i = 1, \ldots, \dim G$,
the left coproduct is defined by $\tilde{\Delta}_* (v_i ) := \tilde{\Delta} v_i + h_i \otimes v_i$,
$\tilde{\Delta}_* h_i= v_i \otimes h_i - h_i \otimes h_i $ and the right coproduct is defined by
$\Delta_* v_i = \Delta v_i + v_i \otimes h_i$, $\Delta_* h_i= h_i  \otimes v_i - h_i \otimes h_i$.
The finite Markov 
$L$-coalgebra $G_*$ is called {\it{the companion}} of $G$. It is defined
up to an isomorphism.
The directed graph so obtained is called the {\it{companion graph}}.
\begin{center}
\includegraphics*[width=5cm]{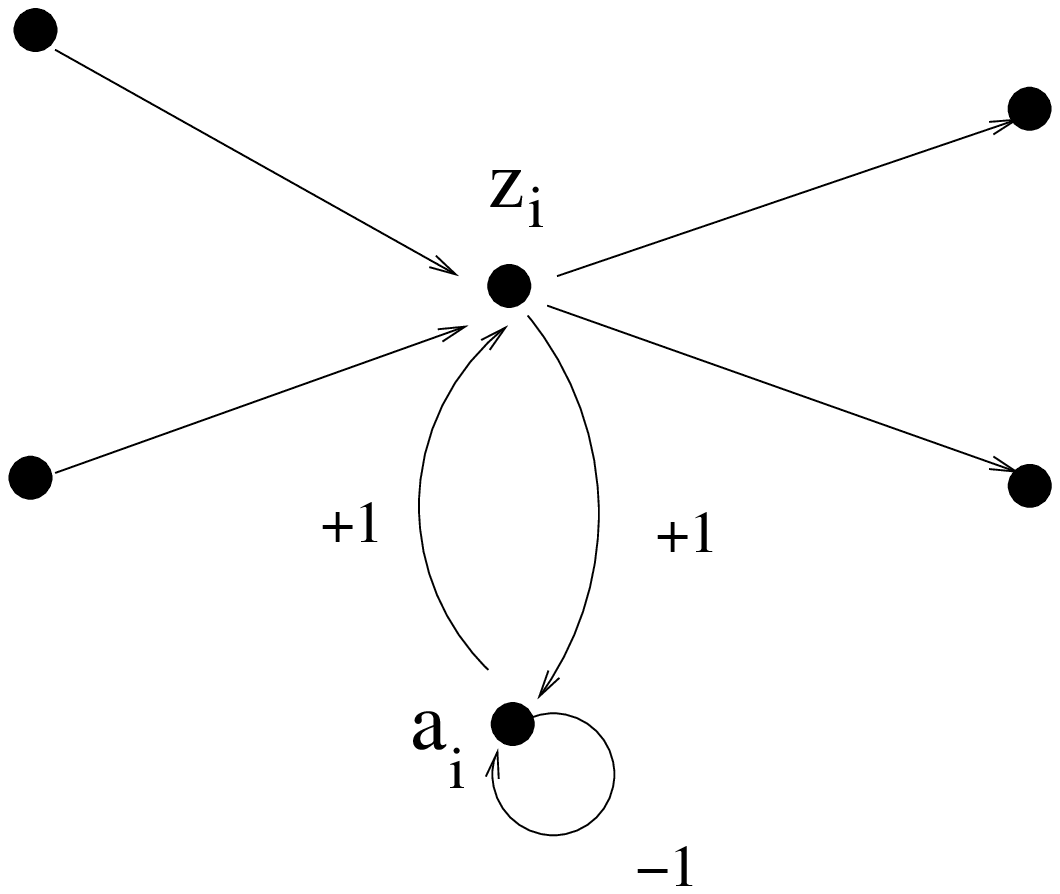}

\textbf{Geometric representation of the companion graph at $v_i$.}
\end{center}
\end{defi}
\begin{lemm}
\label{lmpri}
Let $G$ be a finite Markov $L$-coalgebra with basis $(v_i)_{1 \leq i \leq \dim G}$ and such that the left coproduct is labelled on finite sets 
$J_{v_i}$.
Denote  by $\tilde{w}$ the family of weights $(\tilde{w}_{v_i})_{{v_i} \in G}$
necessary to the definition of the left coproduct of $G$.
The linear mapping $\Psi_{\tilde{w}} : G_*  \xrightarrow{} G_* $
defined by $v_i \mapsto \sum_{k: \ v_{k} \otimes v_i \in J_{v_i}} \tilde{w}_{v_{i}}(v_{k} \otimes v_i) \ v_{k} +h_i$, $h_i \mapsto \sum_{k: \ v_{k} \otimes v_i \in J_{v_i}} \tilde{w}_{v_{i}}(v_{k}\otimes v_i) \ v_{k} +h_i +v_i$
and the linear mapping $\Phi_{\tilde{w}}: G_*  \xrightarrow{} G_* $ defined by $v_i \mapsto h_i - v_i$ and
$h_i \mapsto \sum_{k: \ v_{k} \otimes v_i \in J_{v_i}} \tilde{w}_{v_{i}}(v_{k}\otimes v_i) \ (v_{k} - h_{k}) + v_i$, for all $i$
such that $1 \leq i \leq \dim G$. Then $\Phi_{\tilde{w}} = (\Psi_{\tilde{w}})^{-1}$.
\end{lemm}
\Proof
Fix $i$ such that $1 \leq i \leq \dim G$
and let us prove that $\Psi_{\tilde{w}}$ is invertible.
\begin{eqnarray*}
v_i &\xrightarrow{\Psi_{\tilde{w}}}& \sum_{k: \ v_{k}\otimes v_i \in J_{v_i}} \tilde{w}_{v_{i}}(v_{k}\otimes v_i) \ v_{k}  + h_i \\ &\xrightarrow{\Phi_{\tilde{w}}}& \sum_{k: \ v_{k}\otimes v_i \in J_{v_i}} \tilde{w}_{v_{i}}(v_{k}\otimes v_i) \ (h_{k} - v_{k}) +
 \sum_{k: \ v_{k}\otimes v_i \in J_{v_i}}  \tilde{w}_{v_{i}}(v_{k}\otimes v_i) \ (v_{k} - h_{k}) + v_i = v_i. \\
v_i &\xrightarrow{\Phi_{\tilde{w}}}& h_i - v_i \\ &\xrightarrow{\Psi_{\tilde{w}}}& \sum_{k: \ v_{k}\otimes v_i \in J_{v_i}} \tilde{w}_{v_{i}}(v_{k}\otimes v_i) \  v_{k} + h_i + v_i  - \left(\sum_{k: \ v_{k}\otimes v_i \in J_{v_i}}  \tilde{w}_{v_{i}}(v_{k}\otimes v_i) \ v_{k} + h_i \right) = v_i.
\end{eqnarray*}
\begin{eqnarray*}
h_i &\xrightarrow{\Psi_{\tilde{w}}}& \sum_{k: \ v_{k}\otimes v_i \in J_{v_i}} \tilde{w}_{v_{i}}(v_{k}\otimes v_i) \ v_{k}  + h_i + v_i \\ &\xrightarrow{\Phi_{\tilde{w}}}& \sum_{k: \ v_{k}\otimes v_i \in J_{v_i}} \tilde{w}_{v_{i}}(v_{k}\otimes v_i) \ ( h_{k}  - v_{k})
+ \left( \sum_{k: \ v_{k}\otimes v_i \in J_{v_i}}  \tilde{w}_{v_{i}}(v_{k}\otimes v_i) \ (v_{k} - h_{k}) + v_i \right) \\
& & + (h_i - v_i) =h_i.\\
h_i &\xrightarrow{\Phi_{\tilde{w}}}&  \sum_{k: \ v_{k}\otimes v_i \in J_{v_i}} \tilde{w}_{v_{i}}(v_{k}\otimes v_i) \  (v_{k} - h_{k}) + v_i \\
&\xrightarrow{\Psi_{\tilde{w}}}&
\sum_{k: \ v_{k}\otimes v_i \in J_{v_i}} \tilde{w}_{v_{i}}(v_{k}\otimes v_i) \  \left(\sum_{l: \ v_{l}\otimes v_k \in J_{v_k}} \tilde{w}_{v_{k}}(v_{l}\otimes v_{k}) \ v_{l} + h_{k} \right) \\ & & -
\sum_{k: \ v_{k}\otimes v_i \in J_{v_i}} \tilde{w}_{v_{i}}(v_{k}\otimes v_i) \ \left( \sum_{l: \ v_{l}\otimes v_k \in J_{v_k}} \tilde{w}_{v_{k}}(v_{l}\otimes v_{k}) \ v_{l} + h_{k} + v_{k} \right) \\ & & + \sum_{k: \ v_{k}\otimes v_i \in J_{v_i}} \tilde{w}_{v_{i}}(v_{k}\otimes v_i) \ v_{k}+ h_i = h_i.
\end{eqnarray*}
\eproof
\NB
There exists an unique couple of automorphisms $(\Psi_{\tilde{w}}, \Phi_{\tilde{w}})$ from $\Aut(V_*)$ such that for all $v_i \in G$, $\tilde{\Delta}_* v_i = \Psi_{\tilde{w}}(v_i) \otimes v_i$
and for all $h_i \in H$, $\tilde{\Delta}_* h_i = - \Phi_{\tilde{w}}(v_i) \otimes h_i$ and $\Phi_{\tilde{w}} = (\Psi_{\tilde{w}})^{-1}$.
What was did with the left coproduct $\tilde{\Delta}_*$ remains exact with
the right coproduct $\Delta_*$. The equations
remain the same except the labels of the sums $ \{k: \ v_{k}\otimes v_i \in J_{v_i}  \}$ which obviously become $ \{k: \ v_{i}\otimes v_k \in I_{v_i} \}$ and where the family of weights $\tilde{w}$ is removed by
the family of weights $w$.
There always exists an unique couple of automorphisms $(\Psi_w, \Phi_w)$ of $\Aut(V_*)$ such that for all $v_i \in G$, $\Delta_* v_i =   v_i \otimes \Psi_w(v_i)$
and for all $h_i \in H$, $\Delta_* h_i = h_i \otimes - \Phi_w (v_i)$ and $\Phi_w = \Psi_w ^{-1} $.
\begin{theo}
With each non $L$-cocommutative finite Markov $L$-coalgebra are associated at least
two representations of braid groups determined by its
coproducts.
\end{theo}
\Proof
With each finite Markov $L$-coalgebra $G$, representing an unique directed graph $G$ is associated, up to an isomorphism, an unique finite Markov
$L$-coalgebra $G_*$, thus an unique companion graph $G_*$.
Therefore, there exists up to an isomorphism, two different automorphisms $\Psi_w$ and $\Psi_{\tilde{w}}$, since the coalgebra is not $L$-cocommutative, from
$\Aut(V_*)$ coding information contained within coproducts of $G_*$. $\hat{\Psi}_w$ and $\hat{\Psi}_{\tilde{w}}$ are solutions of YBE thanks to the theorem \ref{tttqye}.
\eproof
\begin{coro}{[\textsf{General case }]}
Let $V$ be a coalgebra with dimension $\dim V$ and with coproduct $\Delta$. The coproduct
$\Delta$ generates at least one non trivial representation of braid groups.
\end{coro}
\Proof
Let $V$ be a coalgebra with coproduct $\Delta$ such that the
associated directed graph be by hypothesis row and locally finite.
By using the
Sweedler notation \cite{Sweedler}, $\Delta v :=v_{(1)} \otimes v_{(2)}$, for all $v \in V$,
let us introduce the coproducts $\tilde{\Delta} v := v_{(1)} \otimes v$
and $\Delta v :=v \otimes v_{(2)}$. These coproducts define a finite
Markov $L$-coalgebra. If it is not $L$-cocommutative, there are at least two
different representations of braid groups and at least one otherwise.
\eproof

\noindent
\textbf{Acknowledgments:}
The author wishes to thank Dimitri Petritis for useful discussions and fruitful advice for the
redaction of this paper.

\bibliographystyle{plain}
\bibliography{These}

\end{document}